\documentclass[a4paper,11pt,leqno]{amsart}

\usepackage{amsmath,amsfonts,amssymb}
\usepackage{pifont,graphicx,epsfig}
\usepackage[bf]{caption2}
\usepackage{epstopdf}

\def\C{\hbox{\font\dubl=msbm10 scaled 1100 {\dubl C}}}

\def\R{\hbox{\font\dubl=msbm10 scaled 1100 {\dubl R}}}
\def\sC{\hbox{\font\dubl=msbm10 scaled 900 {\dubl C}}}
\def\sR{\hbox{\font\dubl=msbm10 scaled 900 {\dubl R}}}
\def\N{\hbox{\font\dubl=msbm10 scaled 1100 {\dubl N}}}
\def\Re{{\rm{Re}}\,}
\def\Im{{\rm{Im}}\,}

\usepackage{amssymb, amsmath}

\newtheorem{Theorem}{Theorem}

\newtheorem{Lemma}[Theorem]{Lemma}
\newtheorem{Proposition}[Theorem]{Proposition}

\def\eps{\varepsilon}

\usepackage{hyperref}

\title[On the roots of the equation $\zeta(s)=a$]
{On the roots of the equation $\zeta(s)=a$}

\author[R. Garunk\v stis]
{Ram\=unas Garunk\v stis}
\thanks{The first author is supported by grant No MIP-94 from the
 Research Council of Lithuania}

\author[ J. Steuding]
{J\"orn Steuding}

\date{November 2010}

\begin{document}

\begin{abstract}
Given any complex number $a$, we prove that there
are infinitely many simple roots of the equation
$\zeta(s)=a$ with arbitrarily large imaginary
part. Besides, we give a heuristic interpretation
of a certain regularity of the graph of the curve
$t\mapsto \zeta({1\over 2}+it)$. Moreover, we
show that the curve $\sR\ni t\mapsto
(\zeta({1\over 2}+it),\zeta'({1\over 2}+it))$ is
not dense in $\sC^2$.
\end{abstract}

\maketitle

{\small \noindent {\sc Keywords:} Riemann zeta-function, value-distribution}
\\ \bigskip

\section{Introduction and statement of the main results}

The Riemann zeta-function $\zeta(s)$ is of
special interest in number theory and complex
analysis. The real zeros of $\zeta(s)$ are called
trivial; they are located at $s=-2n, n\in\N$. All
other zeros are called nontrivial; they lie in
the critical strip $0<\Re s<1$ and are known to
be relevant in many questions concerning the
distribution of prime numbers. It is well-known
that there are infinitely many nontrivial zeros.
More precisely, for the number $N(T)$ of
nontrivial zeros $\rho=\beta+i\gamma$ satisfying
$0<\gamma\leq T$ the asymptotical formula
$$
N(T)={T\over 2\pi}\log{T\over 2\pi e}+O(\log T)
$$
holds as $T\to\infty$. Conrey \cite{conrey} has
shown that more than two fifths of the zeros lie
on the critical line $\Re s={1\over 2}$ and are
simple; this result has been slightly improved by Bui, Conrey \& Young \cite{bui1}. The famous yet unproved Riemann
hypothesis states that all nontrivial zeros lie
on the critical line and the simplicity
hypothesis claims that all (or at least almost
all) zeros are simple. In this article we are
concerned with the general value-distribution. A
famous open problem in this direction is the
question whether the values $\zeta({1\over
2}+it)$ for $t\in\R$ are dense in the complex
plane.
\smallskip

The zeta-function has no exceptional values (in
the meaning of Nevanlinna theory) except infinity
as was shown by Ye \cite{ye} (see also
\cite{steud}, Chapter 7). There are remarkable
quantitative results. For example, it was shown
by Bohr \& Jessen \cite{bj} that $\log\zeta(s)$
assumes any complex value infinitely often in any
vertical strip $\sigma_1<\Re s<\sigma_2$
satisfying ${1\over 2}<\sigma_1<\sigma_2<1$, and
that for fixed $a\neq 0$ the number of such roots
of the equation $\zeta(s)=a$ with imaginary part
bounded by $T$ has linear asymptotic growth as
$T\to\infty$. For arbitrary complex $a$ the roots
of
$$
\zeta(s)=a
$$
are called $a$-points and are denoted by
$\rho_a=\beta_a+i\gamma_a$. There is an $a$-point
near any trivial zero $s=-2n$ for sufficiently
large $n$ and apart from these $a$-points there
are only finitely many other $a$-points in the
half-plane $\sigma\le0$ (see Lemma \ref{rouch}
below). The $a$-points with $\beta_a\le0$ are
said to be trivial; all other $a$-points are
called nontrivial. For any fixed $a$, there exist
left and right half-planes free of nontrivial
$a$-points (see formulas (\ref{bounded}) and
(\ref{1free})).  As in the case of zeros ($a=0$)
there is a Riemann-von Mangoldt--type formula for
the number $N_a(T)$ of nontrivial $a$-points with
imaginary part $\gamma_a$ satisfying
$0<\gamma_a\leq T$, namely, as $T\to\infty$,
\begin{equation}\label{rvma}
N_a(T)=\frac{T}{2\pi}\log {T\over 2\pi
ec_a}+O(\log T)
\end{equation}
with the constant $c_a=1$ if $a\neq 1$, and
$c_1=2$. This was first proved by Landau
\cite{bll}\footnote{The paper \cite{bll} of Bohr,
Landau \& Littlewood consists of three
independent chapters, the first belonging
essentially to Bohr, the second to Landau, and
the third to Littlewood.} (see also \cite{steud},
Chapter 7). We observe that these asymptotics are
essentially independent of $a$:
$$
N_a(T)\sim N(T).
$$
Levinson \cite{levi} proved that all but
$O(N(T)/\log\log T)$ of the $a$-points with
imaginary part in $T<t<2T$ lie in
\begin{equation}\label{levclus}
\vert \Re s-{\textstyle{1\over
2}}\vert<{(\log\log T)^2\over\log T},
\end{equation}
and hence the $a$-points are clustered around the
critical line. In the special case of zeros this
result was obtained by several mathematicians in
the beginning of the 20th century and was
sometimes misleadingly interpreted as indicator
for the truth of the Riemann hypothesis.
Levinson's result was first proved by Landau
\cite{bll} under assumption of the truth of the
Riemann hypothesis.
\par

We turn to the value-distribution on the critical
line.  It is known (see Corollary 3 in Spira
\cite{spira}) that   $\zeta'(\frac12+it)\ne0$ if
$\zeta(\frac12+it)\ne0.$ Hence, there
are no multiple $a$-points of $\zeta(s)$ on the
critical line $\Re s={1\over 2}$ except for
possibly $a=0$.  It follows that not all combinations
of values for the zeta-function and its
derivative are possible. In this direction the
following theorem is true.

\begin{Theorem}\label{null}
The set
$$\{(\zeta({\textstyle{1\over
2}}+it),\zeta'({\textstyle{1\over
2}}+it))\,:\,t\in\R\}$$ is not dense in $\C^2$.
\end{Theorem}

\noindent It follows that the only possible
singularities of the curve $t\mapsto\zeta({1\over
2}+it)$ have to lie in the origin. By this result
we see that the curve $\R\ni
t\to(\zeta({\textstyle{1\over
2}}+it),\zeta'({\textstyle{1\over 2}}+it))$ fails
to visit all neighborhoods of all points in
$\C^2$.  If the Riemann hypothesis is true, the
values $\zeta(\sigma+it)$ for $t\in\R$ are not
dense for any fixed $\sigma<{1\over 2}$ (see
Proposition \ref{prop} below). As mentioned
above, it is unknown whether the values of the
zeta-function on the critical line are dense in
the complex plane or not. It was shown by Voronin
\cite{vor} that the multidimensional analogue is
true for vertical lines in the open right half of
the critical strip: the set
$\{(\zeta(\sigma+it),\zeta'(\sigma+it),\ldots,\zeta^{(n-1)}(\sigma+it))\,:\,t\in\R\}$
is dense in $\C^n$ for all positive integers $n$
for every fixed $\sigma\in({\textstyle{1\over
2}},1)$. Actually, this result was proved by
Voronin previous to his famous universality
theorem which may be interpreted as an infinite
analogue (see \cite{ste,steud}); the case $n=1$
is due to Bohr \& Courant \cite{boco} (see Figure
1). However, the situation on the critical line
is completely different as follows from
Theorem~\ref{null} above; in particular we see
that Voronin's universality theorem cannot be
extended to any region that covers the critical line.

\begin{figure}[ht]
\includegraphics{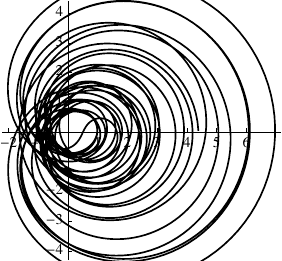}
\hfill
\includegraphics{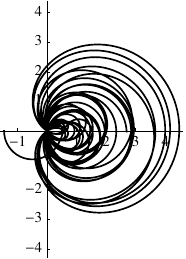}
\hfill
\includegraphics{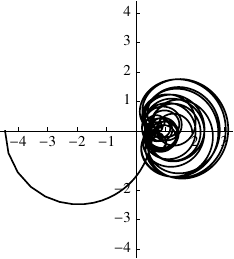}
\caption{\small The curves
$t\mapsto\zeta(\sigma+it)$ for $\sigma={1\over
5},{1\over 2}$, and ${4\over 5}$ from left to
right, all for $t\in[0,100]$. The curve on the
right is known to be dense in the complex plane,
the curve on the left is not dense if Riemann's
hypothesis is true, and for the curve in the
middle this question is open.}
\end{figure}

There is another related topic we want to
investigate. The set of values of $\zeta(s)$ has
the cardinality of the continuum. For some values of $a$
the underlying $a$-points might be non-simple, namely if the
derivative vanishes; however, $\zeta'(s)$ has
only countably many zeros. Thus, there are only
countably many numbers $a$ such that among all
$a$-points there is at least one non-simple
$a$-point. We conjecture that for any fixed
complex number $a$, almost all $a$-points are
simple. By a rather  simple method Conrey, Ghosh
\& Gonek \cite{cgg0} proved that there are
infinitely many simple zeros of the
zeta-function; in \cite{cgg} they have shown by
technical refinement that more than ${19\over
27}$ of the nontrivial zeros are simple provided
the Riemann hypothesis for the Riemann
zeta-function and the generalized Lindel\"of
hypothesis for all Dirichlet $L$-functions are
true; the latter condition has been removed by Bui \& Heath-Brown \cite{bui2}. It is our aim to extend their method to
simple $a$-points; however, for the sake of
simplicity and unconditionality, here we are not
concerned with this type of quantitive results.

\begin{Theorem}\label{uno}
Let $a$ be any fixed complex number. As
$T\to\infty$,
\begin{eqnarray*}
\sum_{0<\gamma_a\leq T} \zeta'(\rho_a)&=&({\textstyle{1\over 2}}-a){T\over 2\pi}\left(\log {T\over 2\pi}\right)^2+(c_0-1+2a){T\over 2\pi}\log {T\over 2\pi}\\
&&+(1-c_0-c_0^2-3c_1-2a){T\over 2\pi}+E(T),
\end{eqnarray*}
where the summation is over nontrivial $a$-points
$\rho_a=\beta_a+i\gamma_a$, the numbers $c_n$ are
the Stieltjes constants (defined by
(\ref{stieltjes}) below), and the error term is
$E(T)\ll T\exp(-C(\log T)^{1/2})$ with some
absolute positive constant $C$; if the Riemann
hypothesis is true, then $E(T)\ll
T^{1/2+\epsilon}$. In any case, for any complex
number $a$ there exist infinitely many simple
$a$-points.
\end{Theorem}

\noindent It should be mentioned that Gonek, Lester \& Milinovich \cite{gonekneu} proved, subject to certain hypotheses, that a positive proportion of the $ a$-points of $\zeta(s)$ are simple; they also obtained an unconditional result for the $ a$-points in fixed strips to the right of the critical line.
\par

In view of the asymptotic formula of
Theorem \ref{uno} the value $a={1\over 2}$
appears to be somehow special for the
zeta-function since in this case the main term is
of lower order. It is easy to see that for
$a={1\over 2}$ the second order term does not
vanish. In fact, the Stieltjes constants are the
coefficients of the Laurent series expansion of
zeta at $s=1$,
\begin{equation}\label{stieltjes}
\zeta(s)=\frac1{s-1}+\sum_{n=0}^\infty
(-1)^n\frac{c_n}{n!}(s-1)^n,
\end{equation}
the constant term
$c_0=\lim_{N\to\infty}(\sum_{n=1}^N{1\over
n}-\log N)=0.577\ldots$ is the Euler-Mascheroni
constant (see Ivi\'c \cite{ivic}). We do not know
why the value ${1\over 2}$ is special in this
sense. Nevertheless, the asymptotics from Theorem
\ref{uno} serve very well for a {\it heuristic}
explanation of the very regular behaviour of the
curve $t\mapsto \zeta({1\over 2}+it)$ as we shall
explain now. Based on computations by Haselgrove
\cite{hasel}, Shanks \cite{shanks} observed that
$\zeta({1\over 2}+it)$ approaches its zeros most
of the times from the 3rd or 4th quadrant,
following Gram's law.\footnote{Recently, it was
shown by Trudgian \cite{gram} that Gram's law
fails for a positive proportion. The first
failure appears at $t=282.454\ldots$.} It was
conjectured by Shanks that the values
$\zeta'({1\over 2}+i\gamma)$ are positive real in
the mean, where $\gamma$ runs through the set of
positive ordinates of the nontrivial zeros. This
follows from the asymptotics obtained by Conrey,
Ghosh \& Gonek \cite{cgg0} under assumption of
the Riemann hypothesis. More precise asymptotical
formulas were derived by Fujii \cite{fujii,fujii3} (see
(\ref{futsch}) below).  Theorem \ref{uno} extends
these results to general $a$.
\begin{figure}[ht]
\centering
\includegraphics[scale=1.1]{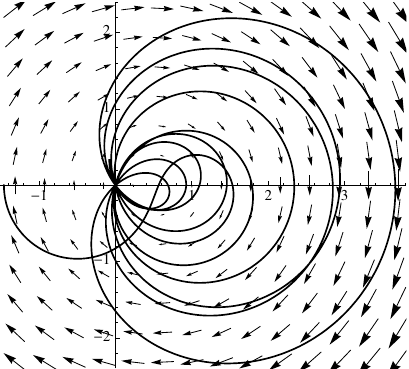}
\caption{\small The curve $t\mapsto \zeta({1\over
2}+it)$ for $0\leq t\leq 50$  and the vector
field $i({1\over 2}-a)$; The value ${1\over 2}$
is the fixed point of the latter -- the eye of
the hurricane!}
\end{figure}
By Levinson's result (\ref{levclus}) almost all
$a$-points lie arbitrarily close to the critical
line, so we may expect that the main contribution
results from these $a$-points. Notice that the
tangent to the curve $t\mapsto\zeta({1\over
2}+it)$ is given by $i\zeta'({1\over 2}+it)$. In
conclusion, the main term $({1\over 2}-a){T\over
2\pi}\log T$ describes how the values
$\zeta({1\over 2}+it)$ approach the value $a$ in
the complex plane on average (see Figure 2).
\par

Finally, we shall prove another theorem of the
same flavour.

\begin{Theorem}\label{due}
For $0\neq \delta:={2\pi\alpha\over \log {T\over
2\pi}}\ll 1$, as $T\to\infty$,
\begin{eqnarray*}
\lefteqn{\sum_{0<\gamma_a\leq T} (\zeta(\rho_a+i\delta)-a)}\\
&=&\left(1-{\sin 2\pi\alpha\over 2\pi\alpha}+i\pi\alpha\left({\sin \pi\alpha\over \pi\alpha}\right)^2-a\{1-\cos 2\pi\alpha+i\sin 2\pi\alpha\}\right){T\over 2\pi}\log {T\over 2\pi}\\
&&+{T\over 2\pi}\Big(-1+\exp(-2\pi i\alpha)\left({1\over i\delta}\left({1\over 1-i\delta}-1\right)-{1\over 1-i\delta}\left(\zeta(1-i\delta)+{1\over i\delta}\right)\right)\\
&&+{\zeta'\over \zeta}(1+i\delta)+{1\over i\delta}-a\Big\{1+\log c_a+(\cos 2\pi\alpha-i\sin 2\pi\alpha)\times\\
&&\times\left({1\over (1-i\delta)^2}-{2\pi
i\alpha\over 1-i\delta}\right)\Big\}\Big)+E(T),
\end{eqnarray*}
uniformly in $\alpha$, where the summation is
over nontrivial $a$-points, $c_a$ is the constant
from (\ref{rvma}) and the error term is of the
same size as in Theorem \ref{uno}.
\end{Theorem}

Theorem \ref{due} generalizes another result of
Fujii \cite{fujii} in the special case $a=0$:
$$
\sum_{0<\gamma\leq T}
\zeta(\rho+i\delta)\sim\left(1-{\sin
2\pi\alpha\over 2\pi\alpha}+i\pi\alpha\left({\sin
\pi\alpha\over \pi\alpha}\right)^2\right){T\over
2\pi}\log T,
$$
where the summation is taken over the zeros
$\rho=\beta+i\gamma$; the precise asymptotical
formula with remainder term is given below as
(\ref{tann}) and (\ref{ca}). A similar discrete
moment was considered by Gonek \cite{gonekk}, who
proved under assumption of the Riemann hypothesis
that
$$
\sum_{0<\gamma\leq
T}\vert\zeta({\textstyle{1\over
2}}+i(\gamma+\delta))\vert^2=\left(1-\left(\displaystyle{\sin(\pi\alpha)\over
\pi\alpha}\right)^2\right){T\over 2\pi}(\log
T)^2+O(T\log T)
$$
uniformly in $\alpha$ for $\vert\alpha\vert \leq
{1\over 2}\log T$. Fujii \cite{fujii2} refined
Gonek's result in replacing the error term by
further explicit main terms plus an error term of
order $O(T^{1/2}(\log T)^3)$. Based on the idea
to model the behaviour of the Riemann
zeta-function on the critical line by the
characteristic polynomials of certain Random
Matrix ensembles, Hughes \cite{hughes}
conjectured, assuming the Riemann hypothesis,
that
$$
\sum_{0<\gamma\leq
T}\vert\zeta({\textstyle{1\over
2}}+i(\gamma+\delta))\vert^{2k} \sim
F_k(2\pi\alpha)a(k)\displaystyle{G(k+1)^2\over
G(2k+1)}{T\over 2\pi}(\log T)^{k^2+1},
$$
where $F_k$ is defined in terms of
Bessel-functions, $a(k)$ is an Euler product, and
$G$ is Barnes' double gamma-function. So far, the
Random Matrix model was not used to do
predictions off the critical line. It would be
very interesting to have a counterpart of Theorem
\ref{due} in Random Matrix Theory.
\par\medskip

The remaining parts of the article are organized
as follows. In the next section we give the proof
of Theorem \ref{null}. In Section 3 we collect
some preliminary results for the proof of Theorem
\ref{uno} which is given in Section 4, resp. the
proof of Theorem \ref{due} which is given in
Section 5. In the final section we state some
concluding remarks. For basic zeta-function
theory we refer to Ivi\'c \cite{ivic}, Titchmarsh
\cite{titch}.

\section{Proof of Theorem \ref{null}}

 Logarithmic
differentiation of the functional equation
\begin{equation}\label{feq}
\zeta(s)=\Delta(s)\zeta(1-s),\qquad
\mbox{where}\quad
\Delta(s)=2^{s}\pi^{s-1}\Gamma(1-s)\sin({\textstyle{\pi
s\over 2}}),
\end{equation}
yields
\begin{equation}\label{01}
\frac{\zeta'}{\zeta}(s)={\Delta'\over\Delta}(s)-\frac{\zeta'}{\zeta}(1-s).
\end{equation}
In view of
\begin{equation}\label{delt}
{\Delta'\over \Delta}(\sigma+it)=-\log{\vert
t\vert\over 2\pi}+O(\vert
t\vert^{-1})\qquad\mbox{for}\quad \vert
t\vert\geq 1
\end{equation}
we get
 \begin{eqnarray}\label{zz}
{\zeta'\over \zeta}({\textstyle{1\over
2}}+it)=-\overline{\left({\zeta'\over
\zeta}\right)}({\textstyle{1\over
2}}+it)-\log{|t|\over 2\pi}+O(|t|^{-1}).
\end{eqnarray}
For $a\in\R$ with $a\neq 0$ we assume that there
is $t_a$ such that
\begin{equation}\label{ineq}
\vert\zeta({\textstyle{1\over
2}}+it_a)-a\vert<\epsilon\qquad\mbox{and}\qquad
\vert\zeta'({\textstyle{1\over
2}}+it_a)-a\vert<\epsilon,
\end{equation}
where $0<\epsilon<\vert a\vert$.  Then
$$
{\zeta'\over \zeta}({\textstyle{1\over
2}}+it_a)={a+O(\epsilon)\over
a+O(\epsilon)}=1+O(\epsilon).
$$
If $a$ is sufficiently large  then  $|t_a|\ge1$.
Hence, we deduce from (\ref{zz}) that
$$
2=-\log{|t_a|\over 2\pi}+O(|t_a|^{-1}+\epsilon).
$$
 For sufficiently large $a$
the latter formula is in contradiction with
(\ref{ineq}). This finishes the proof of Theorem
\ref{null}.

\section{Preliminaries}

In the sequel we write the complex variable as
$s=\sigma+it$ with real $\sigma,t$. We start with
the growth of the zeta-function in the left
half-plane.

\begin{Lemma}\label{3}
There is a constant $c>0$ such that, for
$\sigma\le0$ and $\vert t\vert\geq 2$,
\begin{equation*}
|\zeta(\sigma+it)|>\frac{c\vert
t\vert^{1/2-\sigma}}{\log t}.
\end{equation*}
If Riemann's hypothesis is true, then for any
$\eps>0$ and any $\sigma_0<\frac12$ there is
$c=c(\eps, \sigma_0)>0$ such that, for
$\sigma\le\sigma_0<\frac12$ and $\vert t\vert\geq
2$,
\begin{equation*}
|\zeta(\sigma+it)|>c\vert
t\vert^{1/2-\sigma+\eps}.
\end{equation*}
\end{Lemma}

\noindent {\bf Proof.} It is known (see Patterson
\cite{patt}, Exercise 4.6) that
$$
\zeta(1+it)\gg\frac1{\log t}\quad\mbox{for}\quad
|t|\ge2.
$$
Then the first part of the lemma follows by the
functional equation (\ref{feq}) in combination
with Stirling's formula
\begin{equation}\label{stir}
\Gamma(s)=\exp\left(\left(s-\frac12\right)\log
s-s+\frac{\log
2\pi}{2}\right)\left(1+O_{\sigma_1}\left(\frac1
s\right)\right)\quad\mbox{for}\quad
\sigma\ge\sigma_1>0.
\end{equation}
Assuming the Riemann hypothesis, for fixed
$\sigma>{1\over 2}$ we may use the bound
$\zeta(\sigma+it)\gg \vert t\vert^{-\eps}$ from
\cite{titch}, Chapter 14.2. This finishes the
proof of the lemma.
\medskip

We have the following application of the previous
lemma.
\begin{Proposition}\label{prop}
If the Riemann hypothesis is true, then the
values $\zeta(\sigma+it)$ for $t\in\R$ are not
dense for any fixed $\sigma<{1\over 2}$.
\end{Proposition}

\noindent {\bf Proof.} By Lemma \ref{3}, for
$|t|>T_0>2$,  the values of $|\zeta(\sigma+it)|$
 are greater than some constant
$C=C(T_0)$. Then the curve $\zeta(\sigma+it)$,
$t\in[-T_0,T_0]$, which is continuous and of
finite length, can not be dense in the disc
$|z|\le C$.
\medskip

The next lemma shows that certain $a$-points are
related to trivial zeros of the zeta-function:

\begin{Lemma}\label{rouch}
For any complex number $a$ there exists a
positive integer $N$ such that there is a simple
$a$-point of $\zeta(s)$ in a small neighbourhood
around $s=-2n$ for all positive integers $n\geq
N$; apart from these there are no other
$a$-points in the left half-plane $\Re s\le0$
except possibly finitely many near $s=0$.
\end{Lemma}

\noindent This observation is due to Levinson
\cite{levi}; for the proof one applies the
functional equation for $\zeta$ in combination
with Rouch\'e's theorem and Stirling's formula.
The second assertion follows from Lemma \ref{3}.
\par

Now we investigate the order of growth of the
almost entire function $\zeta(s)-a$. Note that
the order of an entire function $f$ is defined to
be the infimum of all real numbers $b$ for which
the estimate
$$
|f(s)|\le\exp(|s|^b)
$$
holds for all sufficiently large $\vert s\vert$.
The following lemma is well-known in the case
$a=0$; the general case can be treated similarly.

\begin{Lemma}
For any $a$ the function
$(s-1)\left(\zeta(s)-a\right)$ is entire of order
$1$.
\end{Lemma}

\noindent {\bf Proof} is analogous to the proof
of Theorem 2.12 in Titchmarsh \cite{titch}.

\medskip

The next lemma generalizes the well-known partial
fraction decomposition of the logarithmic
derivative of $\zeta$:

\begin{Lemma}\label{5} Let $a$ be a fixed complex number. Then, for $-1\leq\sigma\leq2, |t|\geq 1$,
$$
\frac
{\zeta'(s)}{\zeta(s)-a}=\sum_{|t-\gamma_a|\leq1}{1\over
s-\rho_a}+O(\log(|t|+1)),
$$
where the summation is taken over all $a$-points
$\rho_a=\beta_a+i\gamma_a$ satisfying $\vert
t-\gamma_a\vert\leq 1$.
\end{Lemma}

\noindent {\bf Proof.} Since $(s-1)(\zeta(s)-a)$
is an integral function of order one (by the
previous lemma), Hadamard's factorization theorem
yields
$$
(s-1)(\zeta(s)-a)=\exp(A+Bs)\prod_{\rho_a}\left(1-\frac
s\rho_a\right)\exp\left({\frac s\rho_a}\right),
$$
where $A$ and $B$ are certain complex constants
and the product is taken over {\it all} zeros
$\rho_a$ of $(s-1)(\zeta(s)-a)$ (trivial and
nontrivial $a$-points). Hence, taking the
logarithmic derivative, we get
$$
\frac
{\zeta'(s)}{\zeta(s)-a}=B-\frac1{s-1}+\sum_{\rho_a}\left(\frac1{s-\rho_a}+\frac1{\rho_a}\right);
$$
the latter formula can be found in \cite{bll},
however, for our reasoning we prefer to work with
a truncated version. By Lemma \ref{rouch} there
exists a positive constant $c$ such that the
imaginary parts of all $a$-points in the left
half-plane lie in the interval $[-c,c]$.
Moreover, it follows that there are ${1\over
2}\sigma+O(1)$ many of these trivial $a$-points
with real part greater than $-\sigma$ as
$\sigma\to+\infty$. Thus, for $s$ distant from
any of these trivial $a$-points, we have
\begin{align*}
\sum_{\text{trivial} \
\rho_a}\left(\frac1{s-\rho_a}+\frac1{\rho_a}\right)&\ll\sum_{\text{trivial}
\
\rho_a}\frac{\sqrt{\sigma^2+t^2}}{\sqrt{\beta_a^2+\gamma_a^2}\sqrt{(\sigma-\beta_a)^2+(t-\gamma_a)^2}}\\&\ll1+
\int_1^\infty \frac{|t|}{x\sqrt{x^2+t^2}} d
x\ll\log t
\end{align*}
as $t\to\infty$. Hence, for those values of $s$,
$$
\frac
{\zeta'(s)}{\zeta(s)-a}=B-\frac1{s-1}+\sum_{\text{nontrivial}
\
\rho_a}\left(\frac1{s-\rho_a}+\frac1{\rho_a}\right)+O(\log
t).
$$
Note that the main term in the Riemann-von
Mangoldt type formula for the number of
$a$-points (\ref{rvma}) does not depend on $a$.
Therefore, the same reasoning as for $a=0$, the
case of zeros of $\zeta$, can be applied to the
latter formula (see Titchmarsh \cite{titch}, \S
9.6). This yields the assertion of the lemma.
\medskip

\section{Proof of Theorem \ref{uno}}

By the calculus of residues,
\begin{equation}\label{contour}
\sum_{0<\gamma_a\leq T}\zeta'(\rho_a)={1\over
2\pi i}\oint {\zeta'(s)^2\over \zeta(s)-a} d s,
\end{equation}
where the integration is taken over a rectangular
contour in counterclockwise direction according
to the location of the nontrivial $a$-points of
$\zeta(s)$, to be specified below. In view of the
Riemann-von Mangoldt-type formula (\ref{rvma})
the ordinates of the $a$-points cannot lie too
dense. For any large $T_0$ we can find a
$T\in[T_0,T_0+1)$ such that
\begin{equation}\label{condi}
\min_{\rho_a}\vert T-\gamma_a\vert\geq {1\over
\log T},
\end{equation}
where the minimum is taken over all nontrivial
$a$-points $\rho_a=\beta_a+i\gamma_a$. We shall
distinguish the cases $a\neq 1$ and $a=1$.
\par\medskip

First, lets assume that $a\neq 1$. We choose
$B=\log T$. For $\sigma\to+\infty$ we have that
\begin{equation}\label{bounded}
\zeta(\sigma+it)=1+o(1)
\end{equation}
uniformly in $t$. Thus there are no $a$-points in
the half-plane $\Re s>B-1$. Further, define
$b=1+{1\over \log T}$. Then we may suppose that
there are no $a$-points on the line segments
$[B,B+iT]$ and $[1-b,1-b+iT]$ (by varying $b$
slightly if necessary). Moreover, in view of
Lemma \ref{rouch} there are only finitely many
trivial $a$-points to the right of $\Re s=1-b$.
Hence, in (\ref{contour}) we may choose the
counterclockwise oriented rectangular contour
${\mathcal R}$ with vertices $1-b+i,B+i, B+iT,
1-b+iT$ at the expense of a small error for
disregarding the finitely many nontrivial
$a$-points below $\Im s=1$ and for counting
finitely many trivial $a$-points to the right of
$\Re s=1-b$ :
$$
\sum_{0<\gamma_a\leq T}\zeta'(\rho_a)={1\over
2\pi i}\int_{\mathcal R} {\zeta'(s)^2\over
\zeta(s)-a} d s+O(1).
$$
If there is any $a$-point on the line segment
$[1-b+i,B+i]$, we exclude this value by a small
indention; the contribution of the integral over
this interval is bounded, hence negligible.
\par

Next we consider the integral over the upper
horizontal line segment $[B+iT,1-b+iT]$. By the
Phragm\'en-Lindel\"of principle and by the
functional equation (\ref{feq}), for $
\sigma\ge-3$,
\begin{equation}\label{phrag}
\zeta(\sigma+it)\ll \vert
t\vert^{\max\{(1-\sigma),
0\}/2+\epsilon}\qquad\mbox{as}\quad \vert
t\vert\to\infty
\end{equation}
with an implicit constant depending only on
$\eps$ (see Titchmarsh \cite{titch}, \S 5.1).
Hence by Cauchy's integral formula we deduce, for
$ \sigma\ge-2$,
\begin{equation*}
\zeta'(\sigma+it)\ll \vert
t\vert^{\max\{(1-\sigma),
0\}/2+\epsilon}\qquad\mbox{as}\quad \vert
t\vert\to\infty.
\end{equation*}
Then from Lemma \ref{5} in view of the number of
nontrivial $a$-points (\ref{rvma}) we get, for
$\sigma\ge1-b$,
\begin{equation}\label{00}
{\zeta'(\sigma+iT)^2\over \zeta(\sigma+iT)-a}\ll
T^{(1-\sigma)/2+\epsilon}.
\end{equation}
Consequently, the integrals over the horizontal
line segments contribute at most
$O(T^{1/2+\epsilon})$.

It remains to consider the vertical integrals.
For $\sigma\to+\infty$,
$$\zeta'(s)\ll2^{-\sigma}$$
uniformly in $t$. This and the formula
(\ref{bounded}) give
$$
\int_{B}^{B+iT}{\zeta'(s)^2\over \zeta(s)-a}d
s\ll T^{-2\log 2}\log T.
$$

Collecting together,
\begin{equation}\label{contour2}
\sum_{0<\gamma_a\leq T}\zeta'(\rho_a)=-{1\over
2\pi i}\int_{1-b}^{1-b+iT}{\zeta'(s)^2\over
\zeta(s)-a} d s+O(T^{1/2+\epsilon}).
\end{equation}
Hence, it remains to evaluate the integral over
the left vertical line segment $[1-b+iT,1-b]$ of
${\mathcal R}$.

By Lemma \ref{3} there exists a positive constant
$A$, depending only on $a$, such that
$$
\left\vert{a\over
\zeta(s)}\right\vert<{\textstyle{1\over
2}}\qquad\mbox{for}\quad s=1-b+it,\ \vert
t\vert\geq A.
$$
Hence, the geometric series expansion
$$
{1\over \zeta(s)-a}={1\over
\zeta(s)}\left(1+{a\over
\zeta(s)}+\sum_{k=2}^\infty \left({a\over
\zeta(s)}\right)^k\right)
$$
is valid for $s$ from $[1-b+iA,1-b+iT]$. Since
$$
{1\over 2\pi
i}\int_{1-b+i}^{1-b+iA}{\zeta'(s)^2\over
\zeta(s)-a} d s\ll 1,
$$
we deduce
\begin{eqnarray}\label{inte}
\lefteqn{-{1\over 2\pi i}\int_{1-b}^{1-b+iT}{\zeta'(s)^2\over \zeta(s)-a} d s}\nonumber\\
&=&-{1\over 2\pi i}\int_{1-b+iA}^{1-b+iT}\left\{{\zeta'^2\over \zeta}(s)+a\left({\zeta'\over \zeta}(s)\right)^2+\frac{\zeta'^2}{\zeta}(s)\sum_{k=2}^\infty\left(\frac a{\zeta(s)}\right)^k\right\} d s\nonumber\\
&&+O(1)\\
&=&{\mathcal J}_1+{\mathcal J}_2+{\mathcal
J}_3+O(1),\nonumber
\end{eqnarray}
say. To estimate the third integral, we use Lemma
\ref{3} in combination with (\ref{01}) and
(\ref{delt}) in order to obtain
\begin{eqnarray}\label{j3}
{\mathcal J}_3&=&-{a\over 2\pi i} \int_{1-b+iA}^{1-b+iT}\left({\zeta'\over \zeta}(s)\right)^2\sum_{\ell=1}^\infty \left({a\over \zeta(s)}\right)^\ell d s\nonumber\\
&\ll& T(\log T)^2\sum_{\ell=1}^\infty \left({\log
T\over T^{1/2}}\right)^\ell\ll T^{1/2}(\log T)^3.
\end{eqnarray}

Applying the functional equation in the form
(\ref{01}), we find for the second integral in
(\ref{inte})
\begin{eqnarray*}
{\mathcal J}_2&=&-{a\over 2\pi i}\int_{1-b+iA}^{1-b+iT}\left({\Delta'\over \Delta}(s)-{\zeta'\over \zeta}(1-s)\right)^2 d s\\
&=&-{a\over 2\pi i}\int_{1-b+iA}^{1-b+iT}\left\{\left({\Delta'\over \Delta}(s)\right)^2-2{\Delta'\over \Delta}(s){\zeta'\over \zeta}(1-s)+\left({\zeta'\over \zeta}(1-s)\right)^2\right\} d s\\
&=&{\mathcal H}_1+{\mathcal H}_2+{\mathcal H}_3,
\end{eqnarray*}
say. In combination with the asymptotic formula
(\ref{delt}) we easily get
\begin{eqnarray*}
{\mathcal H}_1&=&-{a\over 2\pi}\int_{A}^{T}\left(-\log {t\over 2\pi}+O(t^{-1})\right)^2 d t\\
&=&-a\left\{{T\over 2\pi}\left(\log {T\over
2\pi}\right)^2-{T\over \pi}\log {T\over
2\pi}+{T\over \pi}\right\}+O(\log T).
\end{eqnarray*}
In a similar way we find
\begin{eqnarray*}
{\mathcal H}_2&=&-{a\over \pi}\int_{A}^{T}\left(\log {t\over 2\pi}+O(t^{-1})\right){\zeta'\over \zeta}(b-it) d t\\
&=&{a\over \pi}\sum_{m=2}^\infty {\Lambda(m)\over
m^b}\int_A^{T}\left(\log {t\over
2\pi}+O(t^{-1})\right)\exp(-it\log m) d t.
\end{eqnarray*}
Integration by parts shows that the integral is
$O(\log T)$, hence
$$
{\mathcal H}_2\ll \log
T\sum_{m=2}^\infty{\Lambda(m)\over m^b}=\log
T\left\vert{\zeta'\over \zeta}(b)\right\vert\ll
(\log T)^2.
$$
The same reasoning shows
$$
{\mathcal H}_3\ll \sum_{m,n=2}^\infty
{\Lambda(m)\Lambda(n)\over
(mn)^b}\left\vert\int_{A}^{T}\exp(-it\log(mn)) d
t\right\vert\ll (\log T)^2.
$$
Hence, collecting together we find
\begin{eqnarray}\label{j2}
{\mathcal J}_2&=&-a\left\{{T\over
2\pi}\left(\log{T\over 2\pi}\right)^2-{T\over
\pi}\log{T\over 2\pi}+{T\over
\pi}\right\}+O((\log T)^2).
\end{eqnarray}

It remains to evaluate the first integral on the
right-hand side of (\ref{inte}). It should be
noted that this is essentially the integral
giving the main term in the proof of Conrey,
Ghosh \& Gonek for the existence of infinitely
many simple zeros \cite{cgg0}, resp. \cite{cgg}
(apart from the mollifier used there to obtain
conditional quantitive results), so
$$
{\mathcal J}_1=-{1\over 2\pi
i}\int_{b}^{b+iT}{\zeta'^2\over \zeta}(s) d
s=\sum_{0<\gamma\leq
T}\zeta'(\rho)+O(T^{1/2+\epsilon})\sim {T\over
4\pi}(\log T)^2,
$$
where the summation is over all nontrivial zeros.
Fujii \cite{fujii, fujii3} obtained a more precise
asymptotic formula, namely, as $T\to\infty$,
\begin{equation}\label{futsch}
\sum_{0<\gamma\leq T}\zeta'(\rho)={T\over
4\pi}\left(\log {T\over
2\pi}\right)^2+(c_0-1){T\over 2\pi}\log{T\over
2\pi}+(1-c_0-c_0^2-3c_1){T\over 2\pi}+{\mathcal E}(T),
\end{equation}
where the numbers $c_0,c_1$ are the Stieltjes
constants given by (\ref{stieltjes}), and the
error term ${\mathcal E}(T)$ is $O(T\exp(-C(\log
T)^{1/2}))$ unconditionally with some absolute
positive constant $C$, resp. $O(T^{1/2}(\log
T)^{7/2})$ under assumption of the Riemann
hypothesis; note that \cite{fujii3} contains a correction of the corresponding formula in \cite{fujii}.
\par

Substituting (\ref{j3}), (\ref{j2}), and
(\ref{futsch}) into (\ref{inte}), leads via
(\ref{contour2}) to the desired asymptotical
formula for every $T$ satisfying condition
(\ref{condi}). To get this uniformly in $T$ we
allow an arbitrarily $T$ at the expense of an
error $\ll T^{1/2+\epsilon}$ (by shifting the
path of integration using (\ref{phrag})). This
proves Theorem \ref{uno} in the case $a\neq 1$.
\par\medskip

For $a=1$ we consider the function
$f(s):=2^s(\zeta(s)-1)$ in place of $\zeta(s)-a$.
By the Dirichlet series expansion
\begin{equation}\label{1free}
2^s(\zeta(s)-1)=1+\sum_{n=3}^\infty\left({2\over
n}\right)^s
\end{equation}
it follows that there is a zero-free right
half-plane for $f(s)$. Computing the logarithmic
derivative,
$$
{f'\over f}(s)=\log 2+{\zeta'(s)\over
\zeta(s)-1},
$$
we observe that the non-constant term corresponds
to the logarithmic derivative in the case $a\neq
1$ while the constant term does not contribute by
integration over a closed contour. This proves
Theorem \ref{uno} for general $a$.

\section{Proof of Theorem \ref{due}}

The proof is rather similar to the previous one.
Let the quantities $b,B$ and $A$ be defined as
above. We start with
\begin{equation}\label{start}
\sum_{0<\gamma_a\leq
T}\zeta(\rho_a+i\delta)={1\over 2\pi i}\oint
{\zeta'(s)\over \zeta(s)-a}\zeta(s+i\delta) d s,
\end{equation}
where the integration is again over a rectangular
contour with vertices $1-b+i,B+i,B+iT,1-b+iT$ in
counterclockwise direction. As before we assume
that there are no $a$-points on this contour,
otherwise we can circumvent these values by a
small indention at the expense of an error
$O(T^{1/2+\epsilon})$. By the same reasoning as
above (see (\ref{inte})) the main contribution
comes from the integral
\begin{equation}\label{tzui}
-{1\over 2\pi
i}\int_{1-b+iA}^{1-b+iT}{\zeta'\over
\zeta}(s)\left(1+{a\over
\zeta(s)}+\sum_{k=2}^\infty \left({a\over
\zeta(s)}\right)^k\right)\zeta(s+i\delta) d s,
\end{equation}
defining a sum consisting of three terms. The
first term was essentially already computed by
Fujii \cite{fujii}, when he proved in the case
$a=0$ the asymptotical formula
\begin{eqnarray}\label{tann}
\sum_{0<\gamma\leq T} \zeta(\rho+i\delta)&=&\left(1-{\sin 2\pi\alpha\over 2\pi\alpha}+i\pi\alpha\left({\sin \pi\alpha\over \pi\alpha}\right)^2\right){T\over 2\pi}\log {T\over 2\pi}\nonumber\\
&&+{T\over 2\pi}c(\alpha,T)+{\mathcal E}(T)
\end{eqnarray}
where the summation is taken over the zeros
$\rho=\beta+i\gamma$,
\begin{eqnarray}\label{ca}
c(\alpha,T)&:=&-1+\exp(-2\pi i\alpha)\left\{{1\over i\delta}\left({1\over 1-i\delta}-1\right)-{1\over 1-i\delta}\left(\zeta(1-i\delta)+{1\over i\delta}\right)\right\}\nonumber\\
&&+{\zeta'\over \zeta}(1+i\delta)+{1\over
i\delta},
\end{eqnarray}
and the error term estimate ${\mathcal E}(T)\ll
T\exp(-C(\log T)^{1/2})$ unconditionally, resp.
$\ll T^{1/2}(\log T)^{5/2}$ under assumption of
Riemann's hypothesis. This yields
\begin{eqnarray}\label{kabel}
\lefteqn{-{1\over 2\pi i}\int_{1-b+iA}^{1-b+iT}{\zeta'\over \zeta}(s)\zeta(s+i\delta) d s=\sum_{0<\gamma\leq T} \zeta(\rho+i\delta)+O(T^{1/2+\epsilon})}\nonumber\\
&=&\left(1-{\sin 2\pi\alpha\over
2\pi\alpha}+i\pi\alpha\left({\sin \pi\alpha\over
\pi\alpha}\right)^2\right){T\over 2\pi}\log
{T\over 2\pi}+{T\over 2\pi}c(\alpha,T)+{\mathcal
E}(T).
\end{eqnarray}

The third term in (\ref{tzui}) contributes again
to the error term. Applying the functional
equations (\ref{feq}) and (\ref{01}), the second
term can be rewritten as
\begin{eqnarray*}
\lefteqn{-{a\over 2\pi i}\int_{1-b+iA}^{1-b+iT}{\zeta'\over\zeta}(s){\zeta(s+i\delta)\over \zeta(s)} d s}\\
&=&-{a\over 2\pi i}\int_{1-b+iA}^{1-b+iT}\left({\Delta'\over\Delta}(s)-{\zeta'\over \zeta}(1-s)\right){\Delta(s+i\delta)\over \Delta(s)}{\zeta(1-s-i\delta)\over \zeta(1-s)} d s\\
&=&{\mathcal K}_1+{\mathcal K}_2,
\end{eqnarray*}
say. It follows from Stirling's formula
(\ref{stir}) that
$$
{\Delta(\sigma+i(\delta+t))\over
\Delta(\sigma+it)}=\exp\left(-i\delta\log{t\over
2\pi}\right)(1+O(t^{-1})).
$$
 Using (\ref{delt}) we have
that
$$
{\mathcal K}_1={a\over 2\pi}\sum_{m,n=1}^\infty
{\mu(m)n^{i\delta}\over
(mn)^b}\int_{A}^{T}\left(\log {t\over
2\pi}+O(t^{-1})\right)\exp\left(it\log(mn)-i\delta\log{t\over
2\pi}\right) d t.
$$
For $mn\neq 1$ the integral can be estimated by
integrating by parts; these terms contribute an
error term $O((\log T)^3)$. Computing the
integral for $m=n=1$ yields
\begin{eqnarray*}
{\mathcal K}_1&=&{a\over 2\pi}\int_{A}^{T}\log
{t\over 2\pi}\exp\left(-i\delta\log{t\over
2\pi}\right) d t
+O((\log T)^3)\\
&=&a\exp(-2\pi i\alpha){T\over
2\pi}\left(\log{T\over 2\pi}-{1\over
1-i\delta}\right){1\over 1-i\delta}+O((\log
T)^3).
\end{eqnarray*}
This gives besides the first term in (\ref{tzui})
a further contribution to the main term.

Similarly we get ${\mathcal K}_2\ll(\log T)^3$.

Note that
$$
{1\over 1-i\delta}\log{T\over 2\pi}=\log{T\over
2\pi}+{2\pi i\alpha\over 1-i\delta}.
$$
Together with (\ref{kabel}) we get the asymptotic
formula for (\ref{start}). Subtracting
(\ref{rvma}) the proof of Theorem \ref{due} is
complete.
\medskip

It should be noted that differentiation of the
formula of Theorem \ref{due} with respect to
$\alpha$ leads to the formula of Theorem
\ref{uno}; for this purpose one has to be aware
that all error terms are uniform in $\alpha$.

\section{Concluding remarks}

i) Similar graphs as for curves
$t\mapsto\zeta(\sigma+it)$ (Figure 1) appear for
other zeta- and $L$-functions too (see for
example Akiyama \& Tanigawa \cite{akiyama} for
$L$-functions associated with elliptic curves).
It seems that the shape of these curves depends
on the type of functional equation, the location
of zeros, as well as on the first coefficient of
the Dirichlet series expansion.
\par\medskip

ii) It is possible to consider short intervals
$(T,T+H]$ for the imaginary parts of nontrivial
$a$-points in place of $(0,T]$ as was done in
\cite{steu} for zeros; here {\it short} means
that $T^{1/2+\epsilon}\leq H\leq T$. Moreover, we
observe that, generalizing Theorem \ref{uno} and
\ref{due},
$$
\sum_{\vert\gamma_a\vert\leq T}
f(\rho_a)=\sum_{\vert\gamma\vert\leq T}
f(\rho)-{a\over 2\pi i}\int {\zeta'(s)\over
\zeta(s)-a}f(s) d s+\mbox{error},
$$
where the summation on the right-hand side is
taken over nontrivial zeros (and not $a$-points),
and where $f$ is {\it any sufficiently smooth}
Dirichlet polynomial or series. These extensions
will be considered in a sequel to this article.
\par\medskip

iii) As already pointed out in the introduction,
quite much is known about the distribution of
$a$-points to the right of the critical line
whenever $a$ is a fixed complex number different
from zero. On the contrary, on the critical line
for no complex number $a$ apart from zero it is
proved that there exist infinitely many
$a$-points. Hence, it seems to be an interesting
problem to study the location of $a$-points to
the left or to the right of curves
$s=\sigma(t)+it$ with
$\lim_{t\to\infty}\sigma(t)={1\over 2}$. Selberg
\cite{selber} was the first to obtain results on
the statistical distribution of $a$-points in
such regions (even for elements in the Selberg
class).
\par

\par\bigskip

\tiny

\noindent
Ram\= unas Garunk\v stis\\
Faculty of Mathematics and Informatics, Vilnius University\\
Naugarduko 24, 03225 Vilnius, Lithuania\\
ramunas.garunkstis@mif.vu.lt
\smallskip

\noindent
J\"orn Steuding\\
Department of Mathematics, W\"urzburg University\\
Emil-Fischer-Str. 40, 97\,218 W\"urzburg, Germany\\
steuding@mathematik.uni-wuerzburg.de


\medskip

\small

\begin{thebibliography}{9}

\bibitem{akiyama}{\sc S. Akiyama, Y. Tanigawa}, Calculation of values of $L$-functions associated to elliptic curves, {\it Math. Comp.} {\bf 68} (1999), 1201-1231

\bibitem{boco}{\sc H. Bohr, R. Courant}, Neue Anwendungen der Theorie der diophantischen Approximationen auf die Riemannsche Zetafunktion, {\it J. reine angew. Math.} {\bf 144} (1914), 249-274

\bibitem{bj}{H. Bohr, B. Jessen}, \"Uber die Werteverteilung der Riemannschen Zetafunktion, zweite Mitteilung, {\it Acta Math.} {\bf 58} (1932), 1-55

\bibitem{bll}{\sc H. Bohr, E. Landau, J.E. Littlewood},
Sur la fonction $\zeta(s)$ dans le voisinage de la
droite $\sigma = \frac12$, {\it Bull. de l'Acad. royale
de Belgique} (1913), 3-35

\bibitem{bui2}{\sc H.M. Bui, D.R. Heath-Brown}, On simple zeros of the Riemann zeta-function, {\it Bull. Lond. Math. Soc.} {\bf 45} (2013), 953-961

\bibitem{bui1}{\sc H.M. Bui, J.B. Conrey, M.P. Young}, More than 41 \% of the zeros of the zeta function are on the critical line, {\it Acta Arith.} {\bf 150} (2011), 35-64 

\bibitem{conrey}{\sc J.B. Conrey}, More than two fifths of the zeros of the Riemann zeta-function are on the critical line, {\it J. reine angew. Math.} {\bf 399} (1989), 1-26

\bibitem{cgg0}{\sc J.B. Conrey, A. Ghosh, S.M. Gonek}, Simple zeros of zeta functions, {\it Colloq. de Theorie Analyt. des Nombres}, (1985), 77-83

\bibitem{cgg}{\sc J.B. Conrey, A. Ghosh, S.M. Gonek}, Simple zeros of the Riemann zeta-function, {\it Proc. Lond. Math. Soc.}, III. Ser. {\bf 76} (1998), 497-522

\bibitem{fujii}{\sc A. Fujii}, On a conjecture of Shanks, {\it Proc. Japan Acad.} {\bf 70} (1994), 109-114

\bibitem{fujii2}{\sc A. Fujii}, On a mean value theorem in the theory of the Riemann zeta-function, {\it Comment. Math. Univ. St. Paul.} {\bf 44} (1995), 59-67

\bibitem{fujii3}{\sc A. Fujii}, On the distribution of values of the derivative of the Riemann zeta function at its zeros. I, {\it Tr. Mat. Inst. Steklova} {\bf 276} (2012), Teoriya Chisel, Algebra i Analiz, 57-82; translation in {\it Proc. Steklov Inst. Math.} {\bf 276} (2012), 51-76

\bibitem{ste}{\sc R. Garunk\v stis, A. Laurin\v cikas, K. Matsumoto, J. Steuding, R. Steuding}, Effective uniform approximation by the Riemann zeta-function, {\it Publ. Mat.} {\bf 54} (2010), 209-219

\bibitem{gonekk}{\sc S.M. Gonek}, Mean values of the Riemann zeta-function and its derivatives, {\it Inventiones math.} {\bf 75} (1984), 123-141

\bibitem{gonekneu}{\sc S.M. Gonek, S.J. Lester, M.B. Milinovich}, A note on simple $a$-points of $L$-functions, {\it Proc. Amer. Math. Soc.} {\bf 140} (2012), no. 12, 4097-4103

\bibitem{hasel}{\sc C.B. Haselgrove}, {\it Tables of the Riemann Zeta function}, Roy. Soc. Math. Tables. Vol. 6., Cambridge University Press. XXIII (1960)

\bibitem{hughes}{\sc C.P. Hughes}, Random matrix theory and discrete moments of the Riemann zeta-function, {\it J. Phys. A} {\bf 36} (2003), 2907-2917

\bibitem{ivic}{\sc A. Ivi\'c}, {\it The Riemann zeta-function}, John Wiley \& Sons, New York 1985

\bibitem{levi}{\sc N. Levinson}, Almost all roots of $\zeta(s)=a$ are arbitrarily close to $\sigma={\textstyle{1\over 2}}$,
{\it Proc. Nat. Acad. Sci. U.S.A.} {\bf 72}
(1975), 1322-1324

\bibitem{patt}{\sc S.J. Patterson}, {\it An introduction to the theory of  the Riemann zeta-function}, Cambridge Univ. Press, 1995

\bibitem{selber}{\sc A. Selberg}, Old and new conjectures and results about a class of Dirichlet series, in: `{\it Proceedings of the Amalfi Conference on Analytic Number Theory}', Maiori 1989, E. Bombieri et al. (eds.), Universit\`a di Salerno 1992, 367-385

\bibitem{shanks}{\sc D. Shanks}, Review of Haselgrove, {\it Math. Comp.} {\bf 15} (1961), 84-86

\bibitem{spira}{\sc R. Spira}, On the Riemann zeta function, {\it J. London Math. Soc.} {\bf 44} (1969), 325-328

\bibitem{steu}{\sc J. Steuding}, On simple zeros of the
Riemann zeta-function in short intervals on the critical
line, {\it Acta Math. Hung.} {\bf 96} (2002), 259-308

\bibitem{steud}{\sc J. Steuding}, {\it Value-Distribution
of $L$-Functions}, Lecture Notes in Mathematics ,
Vol. 1877, Springer, 2007

\bibitem{titch}{\sc E.C. Titchmarsh}, {\it The theory of the Riemann zeta-function}, 2nd ed., revised by D.R. Heath-Brown, Oxford University Press, 1986

\bibitem{gram}{\sc T. Trudgian}, Gram's Law Fails a Positive Proportion of the Time, available as {\sf arXiv:0811.0883} at {\sf http://front.math.ucdavis.edu/0811.0883}

\bibitem{vor}{\sc S.M. Voronin}, The distribution of the non-zero values of the Riemann zeta-function, {\it Izv. Akad. Nauk Inst. Steklov} {\bf 128} (1972), 131-150 (Russian)

\bibitem{ye}{\sc Z. Ye}, The Nevanlinna functions of the Riemann zeta-function, {\it J. Math. Analysis Appl.} {\bf 233} (1999), 425-435


\end{thebibliography}
\end{document}